\documentclass[]{interact}

\usepackage{breqn}

\usepackage{epstopdf}
\usepackage{subfigure}

\usepackage[numbers,sort&compress,merge]{natbib}
\bibpunct[, ]{(}{)}{,}{n}{,}{,}

\theoremstyle{plain}

\theoremstyle{definition}

\theoremstyle{remark}

\begin{document}

\begin{titlepage}
    \begin{center}
        \vspace*{1cm}
            
        \LARGE
        \textbf{Predicting Failure times for some Unobserved Events with Application to Real-Life Data}
            
            
        \vspace{1.5cm}
            
        \textbf{Mahmoud Mansour}
        
        \Large
        Email: Mahmoud.Mansour@bue.edu.eg

        \vspace{1.5cm}
        \textbf{Mohamed Aboshady}
        
        \Large
        Email: Mohamed.Aboshady@bue.edu.eg
            
            
        \vspace{1.8cm}

        \Large
        Basic Science Department\\
        Faculty of Engineering\\
        The British University in Egypt\\
        El-Sherook City\\
        Egypt\\
        
        \vspace{0.8cm}
        26th April 2022
            
    \vspace{1.5cm}
        \textbf{Corresponding author: Mohamed Aboshady}
    
    \end{center}
\end{titlepage}


\title{Predicting Failure times for some Unobserved Events with Application to Real-Life Data}

\author{\name{Mahmoud Mansour and Mohamed Aboshady\thanks{CONTACT M. Aboshady Email: Mohamed.Aboshady@bue.edu.eg}}\affil{Department of Basic Science, Faculty of Engineering, The British University in Egypt, El Sherook City, Cairo, Egypt}}

\maketitle

\begin{abstract}
This study aims to predict failure times for some units in some lifetime experiments. In some practical situations, the experimenter may not be able to register the failure times of all units during the experiment. Recently, this situation can be described by a new type of censored data called multiply-hybrid censored data. In this paper, the linear failure rate distribution is well-fitted to some real-life data and hence some statistical inference approaches are applied to estimate the distribution parameters. A two-sample prediction approach applied to extrapolate a new sample simulates the observed data for predicting the failure times for the unobserved units.
\end{abstract}


\section{Introduction}

In statistical inference, a process in which data from a randomly collected sample from a population is used to estimate the parameters of a population, observations on all the sample elements are not usually available due to being lost, not recorded, or cost and time constraints. The available data is then called a censored sample. Censoring schemes are classified according to the modality of terminating the experiment, terminating the experiment after a pre-specified time, say $T$, will guarantee the duration of the test but ignore the efficiency, since the number of tested components will be random, and this kind of censoring is denoted as Type-I censoring scheme. On the other hand, if the efficiency level is of greater concern than the duration of the experiment, then the experiment is terminated after testing a pre-determined number of components, say $r$. This scheme is denoted as a Type-II censoring scheme where a specific efficiency is guaranteed and the duration of the test is random. A mixture of these two censoring schemes is commonly used and known as hybrid censoring where it gives a more flexible scheme and could be better managed. In the hybrid censoring schemes, a pre-determined number of the tested components, $r$, and a pre-determined time for ending the experiment, $T$, are specified. Assume we are testing $n$ components, let $X_{i,n}$ be the $i^{th}$ ordered failed component, if the test is terminated at a time $T_1=\min\{X_{r,n},T\}$, then the scheme is called Type-I hybrid censored scheme and if the test is terminated at a time $T_2=\max\{X_{r,n},T\}$, the scheme is called Type-II hybrid censored scheme. Several research articles have studied the estimation of unknown parameters for different distribution functions under the Type-I hybrid censored scheme, see \cite{6,13,16,22,28,29}, and under the Type-II hybrid censored scheme, see \cite{3,5,8,12,15,21}.

Due to the number of failures, say $R_i$, that might occur between any two successive observations without observing the exact failure times for these units, it was necessary to introduce a new scheme to enhance the efficiency of these schemes and consider the effect of these losses. This new scheme was first introduced by Lee and Lee \cite{24} and is known as the multiply Type-II hybrid censoring scheme. Jeon and Kang discussed the estimation for parameters for two different distributions under the multiply Type-II hybrid censoring scheme, see \cite{17,18}.

The exponential distribution, which is a special case of Weibull, Gamma, Gompertz,  exponential-geometric, and Linear Failure rate distributions, is considered one of the main life distributions due to its memory-less and constant failure rate properties. Accordingly, it is important in studying the properties of any lifetime phenomena. In the recent decades, some modifications of exponential distributions were introduced, see \cite{7,27,31}.

Probably and as a result of the simplest form of the linear failure rate  (LFR) distribution which is a generalization of the exponential distribution, it arises in the literature on reliability analysis. This distribution, which belongs to the increasing failure rate class, has two parameters, $\alpha$ and $\beta$, and is a special case of the quadratic failure rate model. Apart from the gamma and Weibull distributions, at $t=0$, the LFR model has a positive failure rate while the gamma and Weibull distribution has a failure rate that is equal to zero. Due to its wide application in human survival data, the LFR distribution gained motivation and many authors have studied its properties, see \cite{9}. In this work, we found that the LFR distribution fitted well with real-life survival data.

The random variable $X$ has LFR distribution with probability density function (PDF) and cumulative distribution function (CDF) defined as
\begin{equation}
f(x)=(\alpha+\beta x)e^{-\left( \alpha x+\frac{\beta}{2}x^{2}\right)
},\ \ \ x>0,\ \alpha,\beta>0. \label{1.1}
\end{equation}
and
\begin{equation}
F(x)=1-e^{-\left( \alpha x+\frac{\beta}{2}x^{2}\right)},\ \ \ x>0,\alpha
,\beta>0. \label{1.2}
\end{equation}

In this paper, it is required to propose the classical and Bayesian estimation procedures for the unknown parameters of LFR distribution under the multiply Type-II hybrid censoring scheme. In the case of the Bayesian procedures, the non-informative priors and the informative priors will be considered. Besides that, a simulation study is performed to assess the quality of the choice of the prior distributions and how this choice affects the accuracy of results. And then the two-sample prediction study is applied to the failure data for two models, clinical and industrial applications. In Section 2, the maximum likelihood estimates (MLEs) of the parameters under consideration are obtained in addition to the corresponding approximate confidence intervals (ACIs). Section 3 is devoted to the Bayesian estimation and the MCMC approach. Good estimators can be used in the Bayesian prediction study which is presented in Section 4. A real data set is analyzed in Section 5 for illustration. The simulation study is presented in Section 6 to assess the quality of obtained estimators. Conclusions are given at the end of this paper.

\section{Maximum-likelihood estimation}

Given the data, the log-likelihood functions are considered to be the basis for finding the parameters estimators. Maximum likelihood estimators have many advantages such as satisfying the invariant properties. Also they are asymptotically unbiased, asymptotically normally distributed and have asymptotically minimum variance, see \cite{4,30} for more information on likelihood theory. Under multiply Type-II hybrid censoring scheme and as shown in Figure 1, one can observe one of the following two types of censored data:

\begin{figure}[h]
	\centering
    		\includegraphics[width=120mm,scale=1]{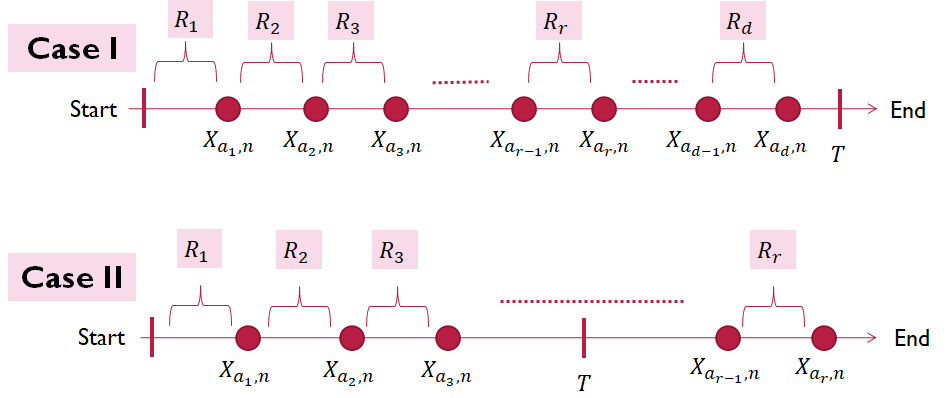}
	\caption{Multiply Type-II hybrid censoring scheme.}
	\label{FIG:1}
\end{figure}

The likelihood function is

\begin{description}
\item Case I: $L\propto\prod\limits_{i=1}^{d}f(x_{a_{i}:n})\prod
\limits_{i=1}^{d-1}\left[  F(x_{a_{i+1}:n})-F(x_{a_{i}:n})\right]  ^{R_{i+1}}\left[  F(x_{a_{1}:n})\right]^{R_{1}}\left[1-F(x_{a_{d::n}})\right]
^{n-d},$

\item Case II: $L\propto\prod\limits_{i=1}^{r}f(x_{a_{i}:n})\prod
\limits_{i=1}^{r-1}\left[F(x_{a_{i+1}:n})-F(x_{a_{i}:n})\right]^{R_{i+1}
}\left[F(x_{a_{1}:n})\right]^{R_{1}}\left[1-F(x_{a_{r::n}})\right]
^{n-r},$
\end{description}

Joining cases I and II, we can rewrite the likelihood function as follows:

\begin{dmath}
L\propto\prod\limits_{i=1}^{m}f(x_{a_{i}:n})\prod\limits_{i=1}^{m-1}\left[F(x_{a_{i+1}:n})-F(x_{a_{i}:n})\right]^{R_{i+1}}\left[F(x_{a_{1}:n})\right]^{R_{1}}\left[1-F(x_{_{a_{m:n}}})\right]^{n-m}, \label{2.1}
\end{dmath}

where $m$ is the number of failure items until the termination point occurred, $R_{i}=a_{i}-a_{i-1}-1$, $a_{0}=0$.\\

The log-likelihood function is
\begin{dmath}
        \ln L\propto\sum\limits_{i=1}^{m}\ln f(x_{a_{i}:n})+\sum\limits_{i=1}^{m-1}R_{i+1}\ln\left[F(x_{a_{i+1}:n})-F(x_{a_{i}:n})\right]+R_{1}\ln\left[F(x_{a_{1}:n})\right]+\left(n-m\right)\ln\left[1-F(x_{_{a_{m:n}}})\right]  \label{2.2}
\end{dmath}

Substituting by the PDF in \eqref{1.1} and the CDF in \eqref{1.2} in the log-likelihood function \eqref{2.2}, we will get:

\begin{dmath}
\ln L\propto\sum\limits_{i=1}^{m}\Bigg[\ln(\alpha+\beta x_{a_{i}:n})-\left(\alpha x_{a_{i}:n}+\frac{\beta}{2}x_{a_{i}:n}^{2}\right)\Bigg]+\sum\limits_{i=1}^{m-1}R_{i+1}\ln\left[e^{-\left(\alpha x_{a_{i}:n}+\frac{\beta}{2}x_{a_{i}:n}^{2}\right)}-e^{-\left(\alpha x_{a_{i+1}:n}+\frac{\beta}{2}x_{a_{i+1}:n}^{2}\right)}\right]+R_{1}\ln\left[1-e^{-\left(\alpha x_{a_{1}:n}+\frac{\beta}{2}x_{a_{1}:n}^{2}\right)}\right]+\left(m-n\right)\left(\alpha x_{_{a_{m:n}}}+\frac{\beta}{2}x_{_{a_{m:n}}}^{2}\right),
\label{2.3}
\end{dmath}

and thus we have the likelihood equations for $\alpha$ and $\beta$
respectively, as

\begin{dmath}
{\sum\limits_{i=1}^{m}}\dfrac{1-x_{a_{i}:n}}{\left[(\alpha+\beta x_{a_{i}:n})-\left(\alpha x_{a_{i}:n}+\frac{\beta}{2}x_{a_{i}:n}^{2}\right)  \right]}+{\sum\limits_{i=1}^{m-1}}\dfrac{R_{i+1}\left[x_{a_{i+1}:n}e^{-\left( \alpha x_{a_{i+1}:n}+\frac{\beta}{2}x_{a_{i+1}:n}^{2}\right)}-x_{a_{i}:n}e^{-\left(\alpha x_{a_{i}:n}+\frac{\beta}{2}x_{a_{i}:n}^{2}\right)}\right]}{\left[e^{-\left(\alpha x_{a_{i}:n}+\frac{\beta}{2}x_{a_{i}:n}^{2}\right)}-e^{-\left(\alpha x_{a_{i+1}:n}+\frac{\beta}{2}x_{a_{i+1}:n}^{2}\right)}\right]}+\dfrac{R_{1}x_{a_{1}:n}\;e^{-\left(\alpha x_{a_{1}:n}+\frac{\beta}{2}x_{a_{1}:n}^{2}\right)}}{\left[1-e^{-\left(\alpha x_{a_{1}:n}+\frac{\beta}{2}x_{a_{1}:n}^{2}\right)}\right]}+\left(m-n\right)
x_{_{a_{m:n}}}=0, \label{2.4}
\end{dmath}

and%

\begin{dmath}
{\sum\limits_{i=1}^{m}}\dfrac{x_{a_{i}:n}\left(1-\frac{1}{2}x_{a_{i}:n}\right)}{\left[(\alpha+\beta x_{a_{i}:n})-\left(\alpha x_{a_{i}:n}+\frac{\beta}{2}x_{a_{i}:n}^{2}\right)\right]}+{\sum\limits_{i=1}^{m-1}}\dfrac{R_{i+1}\left[\frac{1}{2}x_{a_{i+1}:n}^{2}e^{-\left(\alpha
x_{a_{i+1}:n}+\frac{\beta}{2}x_{a_{i+1}:n}^{2}\right)}-\frac{1}{2}x_{a_{i}:n}^{2}e^{-\left(\alpha x_{a_{i}:n}+\frac{\beta}{2}x_{a_{i}:n}^{2}\right)}\right]}{\left[e^{-\left(\alpha x_{a_{i}:n}+\frac{\beta}{2}x_{a_{i}:n}^{2}\right)}-e^{-\left(\alpha x_{a_{i+1}:n}+\frac{\beta}{2}x_{a_{i+1}:n}^{2}\right)}\right]}+\dfrac{\frac{1}{2}R_{1}\;x_{a_{1}:n}^{2}e^{-\left(\alpha x_{a_{1}:n}+\frac{\beta}{2}x_{a_{1}:n}^{2}\right)}}{\left[1-e^{-\left(\alpha
x_{a_{1}:n}+\frac{\beta}{2}x_{a_{1}:n}^{2}\right)}\right]}+\frac{1}{2}\left(m-n\right)x_{a_{m::n}}^{2}=0 \label{2.5}
\end{dmath}

Since it is too difficult to solve the two nonlinear equations \eqref{2.4} and \eqref{2.5} in the two unknown parameters $\alpha$ and $\beta$ simultaneously to get an exact solution, it will be better to solve them numerically using Newton Raphson method to get an approximate solution. For more details on the steps of the Newton Raphson algorithm, see \cite{14}. Finally, the estimates for the parameters $\alpha$ and $\beta$ are the Maximum Likelihood Estimators (MLEs) and will be denoted as $\hat{\alpha}$ and $\hat{\beta}$.

\subsection{Approximate confidence intervals}

The entries of the inverse matrix of the Fisher information matrix $I_{ij}=E\{-[\partial^{2}\ell(\Phi)/\partial\phi_i\partial\phi_j]\}$ gives the asymptotic variances and covariances of the MLEs, $\hat{\alpha}$ and $\hat{\beta}$, where $i,j=1,2$ and $\Phi=(\phi_1,\phi_2)=(\alpha,\beta)$. Unfortunately, it is difficult to obtain the exact closed forms for the previous expectations. Therefore, the observed Fisher information matrix $\hat{I}_{ij}=$ $\left\{  -\left[  \partial^{2}\ell\left(\Phi\right)  /\partial\phi_{i}~\partial\phi_{j}\right]  \right\}  _{\Phi=\hat{\Phi}},~$which is obtained by dropping the expectation operator $E$, will be used to construct confidence intervals for the parameters. The entries of the observed Fisher information matrix are a second partial derivatives of the log-likelihood function, which is easy to be obtained. Hence, the observed information matrix is given by

\begin{equation}
\hat{I}\left(  \alpha,\beta\right)=
\begin{pmatrix}
-\frac{\partial^{2}\ell}{\partial\alpha^{2}} & -\frac{\partial^{2}\ell
}{\partial\alpha\partial\beta}\\
-\frac{\partial^{2}\ell}{\partial\beta\partial\alpha} & -\frac{\partial
^{2}\ell}{\partial\beta^{2}}%
\end{pmatrix}
_{\left(  \alpha,\beta\right)  =\left(  \hat{\alpha},\hat{\beta}\right)}
\label{2.6}
\end{equation}

Therefore, the approximate (or observed) asymptotic variance--covariance matrix $[\hat{V}]$, for the MLEs is obtained by inverting the observed information matrix $\hat{I}(\alpha,\beta)$. Or equivalent

\begin{equation}
[\hat{V}]=\hat{I}^{-1}(\alpha,\beta)=\bigg(\begin{array}[c]{cc}\widehat{Var}(\hat{\alpha}) & cov(\hat{\alpha},\hat{\beta})\\ cov(\hat{\alpha},\hat{\beta}) & \widehat{Var}(\widehat{\beta})\end{array}\bigg)
\label{2.7}
\end{equation}

Under some regularity conditions, it is well known that $(\hat{\alpha},\hat{\beta})$ will be approximately distributed as multivariate normal with mean $(\alpha,\beta)$ and covariance matrix $I^{-1}(\alpha,\beta)$, see \cite{23}.

Then, the $100(1-\gamma)\%$ two sided confidence intervals of $\alpha$ and $\beta,$ can be given by%

\begin{equation}
\widehat{\alpha}\pm Z_{\frac{\gamma}{2}}\sqrt{\widehat{Var}(\hat{\alpha})}\text{ and }\widehat{\beta}\pm Z_{\frac{\gamma}{2}}\sqrt{\widehat{Var}(\widehat{\beta})},\label{2.8}
\end{equation}

where $Z_{\frac{\gamma}{2}}$ is the percentile of the standard normal distribution with right-tail probability $\frac{\gamma}{2}$.

\section{Bayes estimation}

In this section, the Bayesian estimates of the two unknown parameters $\alpha$ and $\beta$ are obtained against the squared error loss function. Assuming that the two parameters $\alpha$ and $\beta$ are independent and follows the jeffrey prior distributions.

\begin{align}
\pi_{1}\left(  \alpha\right)   &  =\alpha^{-1}\ \text{
\ \ \ \ \ \ \ \ \ \ \ \ \ \ \ \ },\ \alpha>0,\nonumber\\
\pi_{2}\left(  \beta\right)   &  =\beta^{-1}\ \text{
\ \ \ \ \ \ \ \ \ \ \ \ \ \ \ \ },\ \beta>0.
\label{3.1}
\end{align}

Most of the researchers in statistical inference use the Bayesian estimation for the unknown parameters as it reduces the posterior expected value for the loss functions, see \cite{1,2,10,25,26}. Using Bayes' theorem and up to proportionality, the posterior distribution of the parameters $\alpha$ and $\beta$, denoted by $\pi^{\ast}(\alpha,\beta\mid\text{data})$, can be obtained by combining the likelihood function \eqref{2.1} with the prior \eqref{3.1} and will be written as

\begin{equation}
\pi^{\ast}(\alpha,\beta\mid\text{data})=\frac{\pi_{1}(\alpha)~\pi_2(\beta)~L(\alpha,\beta\mid\text{data})}{\int\limits_0^{\infty}\int\limits_0^{\infty}\pi_1(\alpha)~\pi_2(\beta)~L(\alpha,\beta\mid\text{data})~d\alpha d\beta}.\label{3.2}
\end{equation}

A square error loss (SEL) function, which is a commonly used function, is a symmetric loss function as it assigns to over and under estimation equal losses. For $\hat{\phi}$ to be the estimator for the parameter $\phi$, then the square error loss function will be defined as 
\begin{equation}
L\left(\phi,\hat{\phi}\right)=\left(\hat{\phi}-\phi\right)^2,\label{3.3}
\end{equation}

Therefore, for any function in $\alpha$ and $\beta$, say $g(\alpha,\beta)$, we can obtain the Bayes estimate under the SEL function as

\begin{equation}
\hat{g}_{BS}\left(  \alpha,\beta\mid\text{data}\right)  =E_{\alpha
,\beta\mid\text{data}}\left[  g(\alpha,\beta)\right]  ,\label{3.4}
\end{equation}

where

\begin{equation}
E_{\alpha,\beta\mid\text{\b{x}}}\left[  g\left(  \alpha,\beta\right)
\right]  =\frac{\int\limits_{0}^{\infty}\int\limits_{0}^{\infty}g\left(
\alpha,\beta\right)  ~\pi_{1}\left(  \alpha\right)  ~\pi_{2}\left(
\beta\right)  ~L(\alpha,\beta\mid\text{data})d\alpha d\beta}
{\int\limits_{0}^{\infty}\int\limits_{0}^{\infty}\pi_{1}\left(  \alpha\right)
~\pi_{2}\left(  \beta\right)  ~L(\alpha,\beta\mid\text{data})~d\alpha d\beta}. \label{3.5}
\end{equation}

As for the complexity to solve the multiple integrals in Equation \eqref{3.5} analytically, it was proposed to use the MCMC approximation method for the samples generation from the joint posterior density function in Equation \eqref{3.2} and then use them for both computing the Bayes estimate of $\alpha$ and $\beta$ and constructing the associated credible intervals. From Equation \eqref{3.2} and up to proportionality, the joint posterior can be written as follow

\begin{dmath}
\pi^{\ast}\left(\alpha,\beta\mid\text{data}\right)\propto\alpha^{-1}\beta^{-1}\prod\limits_{i=1}^{m}(\alpha+\beta x_{a_{i}:n})e^{-\left(\alpha x_{a_{i}:n}+\frac{\beta}{2}x_{a_{i}:n}^{2}\right)}\times\prod\limits_{i=1}^{m-1}\left[e^{-\left(\alpha x_{a_{i}:n}+\frac{\beta}{2}x_{a_{i}:n}^{2}\right)}-e^{-\left(\alpha x_{a_{i+1}:n}+\frac{\beta}{2}x_{a_{i+1}:n}^{2}\right)}\right]^{R_{i+1}}\times\left[1-e^{-\left(\alpha x_{a_{1}:n}+\frac{\beta}{2}x_{a_{1}:n}^{2}\right)}\right]^{R_{1}}\left[e^{(m-n)\left(\alpha x_{a_{m:n}}+\frac{\beta}{2}x_{a_{m:n}}^{2}\right)}\right]\label{3.9}
\end{dmath}

The full conditionals for $\alpha$ and $\beta$ can be written, up to
proportionality, as

\begin{dmath}
\pi_{1}^{\ast}\left(\alpha\mid\beta,\text{data}\right)\propto\alpha^{-1}\prod\limits_{i=1}^{m}(\alpha+\beta x_{a_{i}:n})e^{-\left(\alpha
x_{a_{i}:n}+\frac{\beta}{2}x_{a_{i}:n}^{2}\right)}\prod\limits_{i=1}^{m-1}\left[e^{-\left(\alpha x_{a_{i}:n}+\frac{\beta}{2}x_{a_{i}:n}^{2}\right)}-e^{-\left(\alpha x_{a_{i+1}:n}+\frac{\beta}{2}x_{a_{i+1}:n}^{2}\right)  }\right]^{R_{i+1}}\left[1-e^{-\left(\alpha x_{a_{1}:n}+\frac{\beta}{2}x_{a_{1}:n}^{2}\right)}\right]^{R_{1}}\left[e^{(m-n)\alpha x_{_{a_{m:n}}}}\right]\label{3.10}
\end{dmath}

and

\begin{dmath}
\pi_{2}^{\ast}\left(\beta\mid\alpha,\text{data}\right)\propto\beta^{-1}\prod\limits_{i=1}^{m}(\alpha+\beta x_{a_{i}:n})e^{-\left(\alpha
x_{a_{i}:n}+\frac{\beta}{2}x_{a_{i}:n}^{2}\right)}\prod\limits_{i=1}^{m-1}\left[e^{-\left(\alpha x_{a_{i}:n}+\frac{\beta}{2}x_{a_{i}:n}^{2}\right)}-e^{-\left(\alpha x_{a_{i+1}:n}+\frac{\beta}{2}x_{a_{i+1}:n}^{2}\right)}\right]^{R_{i+1}}\left[1-e^{-\left(\alpha x_{a_{1}:n}+\frac{\beta}{2}x_{a_{1}:n}^{2}\right)}\right]^{R_{1}}\left[e^{\frac{\beta(m-n)}{2}x_{_{a_{m:n}}}^{2}}\right]\label{3.11}
\end{dmath}

It is noted that the conditional posteriors of $\alpha$ and $\beta$ in Equations \eqref{3.10} and \eqref{3.11} are not in standard forms which leads to the difficulty of using Gibbs sampling. Then it will be required to use the Metropolis-Hasting (M-H) sampler for the implementations of MCMC algorithm. A hybrid algorithm with M-H steps for updating $\alpha$ and $\beta$ given the conditional distributions in Equations \eqref{3.10} and \eqref{3.11} is presented below. The steps in this algorithm illustrate the process of the Metropolis-Hastings algorithm within Gibbs sampling:

\begin{enumerate}
    \item Let the initial guess be $\left(\alpha^{(0)},\beta^{(0)}\right)$.
    \item Set $j=1$.
    \item Using the following M-H algorithm, generate $\alpha^{(j)}$ and $\beta^{(j)}$ from $\pi_1^{\ast}\left(\alpha^{(j-1)}\mid\beta^{(j-1)},\text{data}\right)$ and $\pi_2^{\ast}\left(\beta^{(j-1)}\mid\alpha^{(j)},\text{data}\right)$ with the normal proposal distributions $N\left(\alpha^{(j-1)},var\left(\alpha\right)\right)$ and $N\left(\beta^{(j-1)},var\left(\beta\right)\right)$
    \begin{enumerate}
        \item Generate a proposal $\alpha^{\ast}$ from $N\left(\alpha^{(j-1)},var\left(\alpha\right)\right)$ and $\beta^{\ast}$ from $N\left(\beta^{(j-1)},var\left(\beta\right)\right)$.
        \item Evaluate the acceptance probabilities
\begin{align*}
\eta_{\alpha}=\min\left[1,\frac{\pi_{1}^{\ast}\left(\alpha^{\ast}\mid\beta^{(j-1)},\text{data}\right)}{\pi_{1}^{\ast}\left(\alpha^{(j-1)}\mid\beta^{(j-1)},\text{data}\right)}\right], \eta_{\beta}=\min\left[1,\frac{\pi_{2}^{\ast}\left(\beta^{\ast}\mid\alpha^{(j)},\text{data}\right)}{\pi_{2}^{\ast}\left(\beta^{(j-1)}\mid\alpha^{(j)},\text{data}\right)}\right].
\end{align*}
        \item From a Uniform (0,1) distribution, generate $u_1$ and $u_2$.
        \item If $u_1<\eta_{\alpha}$, accept the proposal and set $\alpha^{(j)}=\alpha^{\ast}$, else set $\alpha^{(j)}=\alpha^{(j-1)}$.
        \item If $u_2<\eta_{\beta}$, accept the proposal and set $\beta^{(j)}=\beta^{\ast}$, else set $\beta^{(j)}=\beta^{(j-1)}$.
    \end{enumerate}
    \item Set $j=j+1$.
    \item Repeat Steps $(3)-(5)~N$ times and obtain $\alpha^{(i)}$ and $\beta^{(i)},~i=1,2,...N$.
    \item To compute the credible intervals of $\alpha$ and $\beta$, order $\alpha^{(i)},\beta^{(i)},~i=1,2,...N$ as $\alpha_{(1)}<\alpha_{(2)}<...<\alpha_{(N)}
,~\beta_{(1)}<\beta_{(2)}<...<\beta_{(N)}~$.\\Then the $100(1-\vartheta)\%$ credible intervals of $\varphi=\alpha$ and $\beta$ become $\left(\varphi_{\left(N~\vartheta/2\right)},\varphi_{\left(N~\left(1-\vartheta/2\right)\right)}\right)$.
\end{enumerate}

The choices of the initial values might affect the convergence and for this reason we will discard $M$ simulated points which is known as the burn-in period. The selected samples, $\alpha^{(j)}$ and $\beta^{(j)},~j=M+1,...,N$, for sufficiently large $N$ will then form an approximate posterior samples from which the Bayesian inferences will be developed.

In reference to the SEL function in equation \eqref{3.5}, the approximate Bayes estimates of $\varphi=\alpha,\beta$ will be written as
\begin{equation}
\hat{\varphi}_{BS}=\frac{1}{N-M}\sum\limits_{j=M+1}^{N}\varphi^{(j)},\label{3.12}
\end{equation}

\section{Two sample prediction}

In this section, we will be working on deriving the interval prediction of the future order statistics from a random sample following a LFR distribution based on multiply Type-II hybrid censoring scheme as it may be interested to get the failure times for some observations from a future sample.

Assume that the order statistics from a future random sample of size $m$ be denoted as $Y_{1:m}\leq Y_{2:m}\leq...\leq Y_{m:m}$. Given a probability density function $f(x)$ and a cumulative distribution function $F(x)$ of a continuous distribution. Then for a random sample of size $m$ from this distribution, the marginal density function of the $s^{th}$ order statistic is given by

\begin{dmath}
g_{Y_{s:m}}(y_{s}\mid\boldsymbol{\theta})=\dfrac{m!}{\left(s-1\right)!\left(m-s\right)!}\left[F(y_{s})\right]^{s-1}\left[1-F(y_{s})\right]^{m-s}f(y_{s})=\sum_{q=0}^{m-s}\dfrac{(-1)^{q}\binom{m-s}{q}m!}{\left(s-1\right)!\left(m-s\right)!}\left[F(y_{s})\right]^{s+q-1}f(y_{s}), \label{4.1}
\end{dmath}

where $y_{s}\geq0$ and $\boldsymbol{\theta=}(\alpha,\beta)$, see \cite{11}. Substituting by equations \eqref{1.1} and \eqref{1.2} in \eqref{4.1}, the marginal density function of $Y_{s:m}$ becomes

\begin{dmath}
g_{Y_{s:m}}(y_{s}\mid\alpha,\beta)=(\alpha+\beta y_{s}
)e^{-(\alpha y_{s}+\frac{\beta}{2}y_{s}^{2})}\sum_{q=0}^{m-s}
\dfrac{(-1)^{q}\binom{m-s}{q}m!}{\left(s-1\right)!\left(m-s\right)
!}\left[1-e^{-\left(\alpha y_{s}+\frac{\beta}{2}y_{s}^{2}\right)
}\right]^{s+q-1} \label{4.2}
\end{dmath}

Multiplying the marginal density function in equation \eqref{4.2} and the joint posterior in equation \eqref{3.9}, then integrating this product over the set $\{(\alpha,\beta);~0<\alpha<\infty,~0<\beta<\infty\}$ we can get the predictive posterior density of future observations under Type-II multiply hybrid censoring scheme as follows:

\begin{equation}
g^{\ast}(y_{s}\mid\mathbf{x})=\int_{0}^{\infty}\int_{0}^{\infty}
g_{Y_{s:m}}(y_{s}\mid\mathbf{x})\pi^{\ast}(\alpha,\beta\mid\mathbf{x})d\alpha d\beta . \label{4.3}
\end{equation}

One can note that equation \eqref{4.3} is difficult to be solved analytically. Therefore, we will use MCMC samples through Gibbs algorithm to obtain an estimator for the predictive posterior $g^{\ast}(y_{s}\mid\mathbf{x})$. Assume that $\{(\alpha_i,\beta_i), i = 1, 2,...,N\}$ be the MCMC samples derived from $\pi^{\ast}(\alpha,\beta\mid\mathbf{x})$, then we can get the consistent estimate of $g^{\ast}(y_{s}\mid\mathbf{x})$ as follows\\

\begin{dmath}
g^{\ast}(y_{s}\mid\mathbf{\text{\b{x}}})=\frac{1}{N-M}\sum_{i=M+1}^{N}
(\alpha^{(i)}+\beta^{(i)}y_s)e^{-(\alpha^{(i)}y_{s}+\frac{\beta^{(i)}}{2}y_{s}^{2})}\times\sum_{q=0}^{m-s}
\dfrac{(-1)^{q}\binom{m-s}{q}m!}{(s-1)!(m-s)!}[1-e^{-(\alpha^{(i)}y_{s}+\frac{\beta^{(i)}}{2}y_{s}^{2})}]  ^{s+q-1}. \label{4.4}
\end{dmath}

In addition to the applications that the loss functions have in the Bayesian estimation of parameters, they also have an important role in the Bayesian prediction. We will find, using a loss function, a point predictor that minimizes the risk among all other predictors which is called the Bayesian point predictor. To predict a future observation, we will consider a SEL function and the Bayesian point predictors of $Y_s, 1\leq s\leq N$, denoted as , $\hat{Y}_{SELP}$ under the SEL function to be calculated from

\begin{dmath}
\hat{Y}_{SELP}=\int_{0}^{\infty}y_{s}g^{\ast}(y_{s}|\mathbf{\text{\b{x}}})dy_{s}=\frac{1}{N-M}\sum_{i=M+1}^{N}\int_{0}^{\infty}y_{s}\text{ }g(y_{s}|\alpha^{(i)},\beta^{(i)},\mathbf{\text{\b{x}}})dy_s. \label{4.5}
\end{dmath}

The interval that depends on the information from a previous sample to obtain results for a future sample from a fixed population with a pre-determined probability is called a prediction interval. Given the density function $g(y_s\mid\alpha, \beta,\mathbf{x})$, its distribution function will be given as

\begin{dmath}
G(y_{s}\mid\alpha,\beta,\mathbf{\text{\b{x}}})=\int_{0}^{y_{s}}\sum_{q=0}^{m-s}\dfrac{(-1)^{q}\binom{m-s}{q}m!}{(s-1)!(m-s)!}[F(y_{s})]  ^{s+q-1}f(y_{s})
=\sum_{q=0}^{m-s}\dfrac{(-1)^{q}\binom{m-s}{q}m!}{(s-1)!(m-s)!}\int_{0}^{y_{s}}[F(y_{s})]^{s+q-1}f(y_{s}). \label{4.6}
\end{dmath}

Let $G^{\ast}_{Y_{s:m}}$ denotes the predictive distribution estimator for $y_s$, then a simulation consistent estimator of $G^{\ast}_{Y_{s:m}}(y_s\mid\textbf{x})$ could be written as follows

\begin{dmath}
G^{\ast}(y_{s}\mid\alpha,\beta,\mathbf{\text{\b{x}}})=\frac{1}{N-M}\sum_{i=M+1}^{N}G(y_{s}|\alpha^{(i)},\beta^{(i)},\mathbf{\text{\b{x}}}). \label{4.7}
\end{dmath}

Therefore, the $100(1-\gamma)\%$ confidence interval for the Bayesian predictive interval can be derived after solving the following two equations simultaneously:

\begin{dmath}
G^{\ast}_{Y_{s:m}}(L_{y_{s:m}}|\mathbf{\text{\b{x}}})=1-\frac{\gamma}{2}\text{ and }G^{\ast}_{Y_{s:m}}(U_{y_{s:m}}|\mathbf{\text{\b{x}}})=\frac{\gamma}{2}, \label{4.8}
\end{dmath}

where $L_{y_{s:m}}$ and $U_{y_{s:m}}$ indicate the lower and upper bounds, respectively. It is evident that is not possible to compute the two equations in \eqref{4.8} analytically. Then, the MCMC method is suggested for constructing the Bayesian prediction intervals.

\section{Application to real-life data}

In this section, two real data sets will be discussed to explain the general methods in the previous sections. The first refer to failure times to some industrial data and the second is an application on some clinical data.

\subsection{Example on Industrial data}

First, we consider the classical real data-set in \cite{20} on the times, in operating days, between successive failures of air conditioning equipment in an aircraft. This data is recorded as shown in Table 1. Comparing the empirical distribution for the failure data and the CDF for the LFR distribution as shown in Figure 2, it was found that the Kolmogorov Smirnov (K-S) distance between them is $0.102272$ with $P$-value equals $0.891983$. Therefore, the LFR distribution is fitting well to the mentioned data.\\ 

\begin{tabular}[c]{llllllllll}
\multicolumn{10}{l}{Table 1: Times between failures of air conditioning equipment in an aircraft}\\\hline
$0.417$ & $0.583$ & $0.833$ & $ 0.958$ & $1.000$ & $1.042$ & $ 1.083$ & $1.208$ & $1.833$ & $ 1.833$\\\hline $2.042$ & $2.333$ & $ 2.458$ & $2.500$ & $2.542$ & $ 2.583$ & $2.917$ & $3.167$ & $ 3.292$ & $3.500$\\\hline $3.750$ & $4.208$ & $4.917$ & $5.417$ & $6.500$ & $7.750$ & $8.667$ & $8.667$ & $12.917$\\\hline\\
\end{tabular}

\begin{figure}[h]\centering\includegraphics[width=100mm,scale=1]{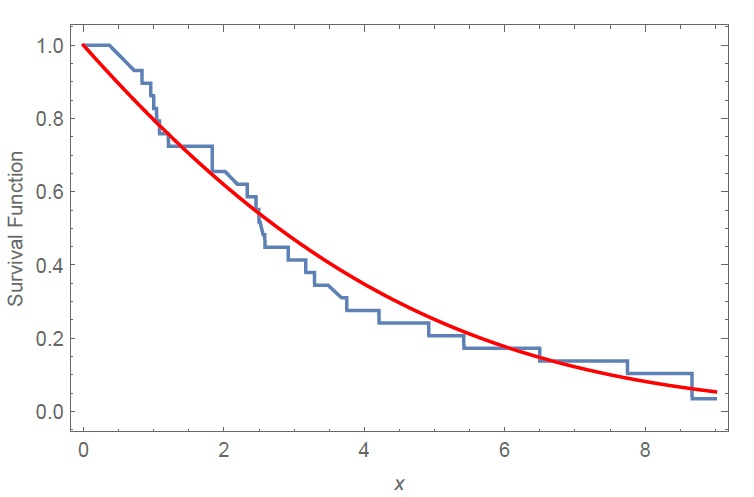}\caption{Empirical and fitted survival functions of data-set in Table 1.}\label{FIG:2}\end{figure}

Then a multiply Type-II hybrid censored sample of effective size $m=23$ was randomly selected from the 29 failure observations in Table 1 with the multiply Type-II hybrid censored scheme $R=(0,1,0,0,1,0,0,0,1,0,0,0,2,0,0,1,0,0,0,0,0,0,0,0,0,0,0,0)$ as shown in Table 2.\\

\begin{tabular}[c]{llllllllll}
\multicolumn{10}{l}{Table 2: Multiply Type-II hybrid failures data}\\\hline
$0.417$ & $0.833$ & $ 0.958$ & $1.042$ & $ 1.083$ & $1.208$ & $ 1.833$ & $2.042$ & $2.333$ & $2.542$\\\hline $2.583$ & $3.167$ & $ 3.292$ & $3.500$ & $3.750$ & $4.208$ & $4.917$ & $5.417$ & $6.500$ & $7.750$\\\hline $8.667$ & $8.667$ & $12.917$\\\hline\\
\end{tabular}

The MLEs for the parameters, $\alpha$ and $\beta$, based on Multiply Type-II hybrid censored data in addition to the Bayes estimates relative to SEL function for the same parameters are displayed in Table 3. Also the $95\%$ approximate confidence intervals and credible intervals for the parameters $\alpha$ and $\beta$ are calculated and the results are also displayed in Table 3.

\begin{center}
\begin{tabular}{llllll}
\multicolumn{6}{l}{Table 3: The point estimates and $95\%$ CIs for $\alpha$ and $\beta$} \\\hline \multicolumn{1}{c}{} & \multicolumn{2}{c}{Point estimate} & \multicolumn{1}{c}{} & \multicolumn{2}{c}{$95\%$ CIs} \\ \cline{2-3} \cline{5-6} \multicolumn{1}{c}{} & \multicolumn{1}{c}{MLE} & \multicolumn{1}{c}{SEL} & \multicolumn{1}{c}{} & \multicolumn{1}{c}{MLE} & \multicolumn{1}{c}{MCMC} \\ \hline $\alpha$ & 0.215785 & 0.214591 & & [0.0317858,0.399783] & [0.212914,0.216269] \\ \hline $\beta$ & 0.0255161 & 0.0246476 & & [-0.0424651,0.0934974] & [0.0244673,0.0248278] \\ \hline
\end{tabular}
\end{center}

From Table 3, the values of estimates are close together which indicates the good performance of the estimators.

Using the MCMC method, a two-sample Bayesian prediction is then computed and its summary statistics is presented in Table 4. It is noted that as $s$ increases, the standard error as well as the confidence intervals widths increase which means that predictive densities get more and more larger as the order statistics become larger.

\begin{center}
\begin{tabular}{llllllllll}
\multicolumn{10}{l}{Table 4: The values of point predictions and $95\%$ PIs for $y_s$.}\\\hline
\multicolumn{1}{c}{} & \multicolumn{1}{c}{} & \multicolumn{3}{c}{95\% PIs} & \multicolumn{1}{c}{} & \multicolumn{1}{c}{} & \multicolumn{3}{c}{95\% PIs} \\ \cline{3-5} \cline{8-10} \multicolumn{1}{c}{$s$} & \multicolumn{1}{c}{SEL} & \multicolumn{1}{c}{Lower} & \multicolumn{1}{c}{Upper} & \multicolumn{1}{c}{Length} & \multicolumn{1}{c}{$s$} & \multicolumn{1}{c}{SEL} & \multicolumn{1}{c}{Lower} & \multicolumn{1}{c}{Upper} & \multicolumn{1}{c}{Length} \\ \hline \multicolumn{1}{c}{1} & \multicolumn{1}{c}{0.251806} & \multicolumn{1}{c}{0.00655} & \multicolumn{1}{c}{0.90774} & \multicolumn{1}{c}{0.90118} & 10 & 3.04697 & 1.5711 & 4.8969 & 3.3258 \\ \hline 2 & 0.51109 & -0.05572 & 1.37661 & 1.43233 & 11 & 3.46596 & 1.86211 & 5.45045 & 3.58833 \\ \hline 3 & 0.778955 & 0.1682 & 1.8007 & 1.6325 & 12 & 3.92753 & 2.18426 & 6.06427 & 3.88001 \\ \hline 4 & 1.01567 & 0.3034 & 2.21117 & 1.90777 & 13 & 4.44512 & 2.5447 & 6.76104 & 4.21634 \\ \hline 5 & 1.34586 & 0.46293 & 2.62059 & 2.15766 & 14 & 5.03974 & 2.95444 & 7.57723 & 4.62279 \\ \hline 6 & 1.64831 & 0.6436 & 3.03674 & 2.39314 & 15 & 5.74726 & 3.43147 & 8.57767 & 5.14621 \\ \hline 7 & 1.96635 & 0.8442 & 3.46591 & 2.62171 & 16 & 6.63705 & 4.00814 & 9.89545 & 5.8873 \\ \hline 8 & 2.30283 & 1.06478 & 3.91421 & 2.84943 & 17 & 7.87495 & 4.75343 & 11.876 & 7.12256 \\ \hline 9 & 2.66143 & 1.3064 & 4.38845 & 3.08205 & 18 & 10.746 & 5.87414 & 15.9654 & 10.0913 \\ \hline
\end{tabular}
\end{center}

\subsection{Example on Clinical data}

The second data-set to be considered here is the survival times (in years) after diagnosis of 43 patients with a certain kind of leukemia as shown in Table 5, see \cite{19}. Comparing the empirical distribution for the failure data and the CDF for the LFR distribution as shown in Figure 3, it was found that the Kolmogorov Smirnov (K-S) distance between them is $0.0816726$ with $P$-value equals $0.914159$. Therefore, the LFR distribution is also fitting well to the this data.\\ 

\begin{tabular}[c]{llllllllll}
\multicolumn{10}{l}{Table 5: Survival times after diagnosis of 43 patients with a kind of leukemia}\\\hline
$0.019$ & $0.129$ & $0.159$ & $0.203$ & $0.485$ & $0.636$ & $0.748$ & $0.781$ & $0.869$ & $1.175$\\\hline $1.206$ & $1.219$ & $1.219$ & $1.282$ & $1.356$ & $1.362$ & $1.458$ & $1.564$ & $1.586$ & $1.592$\\\hline $1.781$ & $1.923$ & $1.959$ & $2.134$ & $2.413$ & $2.466$ & $2.548$ & $2.652$ & $2.951$ & $3.038$\\\hline $3.600$ & $3.655$ & $3.745$ & $4.203$ & $4.690$ & $4.888$ & $5.143$ & $5.167$ & $5.603$ & $5.633$\\\hline $6.192$ & $6.655$ & $6.874$ \\\hline\\
\end{tabular}\\

\begin{figure}[h]
	\centering
		\includegraphics[width=100mm,scale=1]{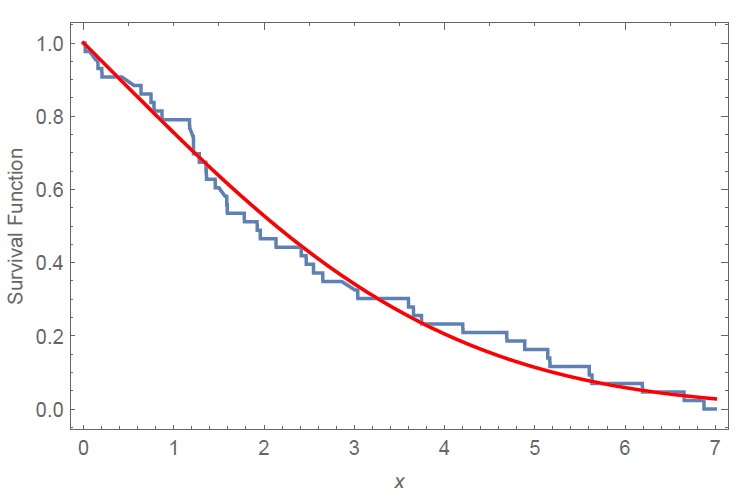}
	\caption{Empirical and fitted survival functions of data-set in Table 5.}
	\label{FIG:3}
\end{figure}

Then a multiply Type-II hybrid censored sample of effective size $m=33$ was randomly selected from the 43 failure observations in Table 5 with the multiply Type-II hybrid censored scheme $R=(0,0,3,0,2,1,1,1,1,0,1,0,0,0,0,0,0,0,0,0,0,0,0,0,0,0,0,0,0,0,0,0,0)$ as shown in Table 6.\\

\begin{tabular}[c]{llllllllll}
\multicolumn{10}{l}{Table 6: Multiply Type-II hybrid failures data}\\\hline
$0.019$ & $0.129$ & $0.636$ & $0.748$ & $1.175$ & $1.219$ & $1.282$ & $1.362$ & $1.564$ & $1.586$\\\hline $1.781$ & $1.923$ & $1.959$ & $2.134$ & $2.413$ & $2.466$ & $2.548$ & $2.652$ & $2.951$ & $3.038$\\\hline $3.600$ & $3.655$ & $3.745$ & $4.203$ & $4.690$ & $4.888$ & $5.143$ & $5.167$ & $5.603$ & $5.633$\\\hline $6.192$ & $6.655$ & $6.874$ \\\hline\\
\end{tabular}\\

The MLEs for the parameters, $\alpha$ and $\beta$, based on Multiply Type-II hybrid censored data in addition to the Bayes estimates relative to SEL function for the same parameters are displayed in Table 7. Also the $95\%$ approximate confidence intervals and credible intervals for the parameters $\alpha$ and $\beta$ are calculated and the results are also displayed in Table 7.

\begin{center}
\begin{tabular}{llllll}
\multicolumn{6}{l}{Table 7: The point estimates and $95\%$ CIs for $\alpha$ and $\beta$} \\\hline \multicolumn{1}{c}{} & \multicolumn{2}{c}{Point estimate} & \multicolumn{1}{c}{} & \multicolumn{2}{c}{$95\%$ CIs} \\ \cline{2-3} \cline{5-6} \multicolumn{1}{c}{} & \multicolumn{1}{c}{MLE} & \multicolumn{1}{c}{SEL} & \multicolumn{1}{c}{} & \multicolumn{1}{c}{MLE} & \multicolumn{1}{c}{MCMC} \\ \hline $\alpha$ & 0.24474 & 0.265739 & & [0.0363186,0.453161] & [0.248641,0.28227] \\ \hline$\beta$ & 0.075861 & 0079844 & & [-0.0617529,0.213475] & [0.0766457,0.0841982] \\ \hline
\end{tabular}
\end{center}

From Table 7, the values of estimates are close together which indicates the good performance of the estimators.

Using the MCMC method, a two-sample Bayesian prediction is then computed and its summary statistics is presented in Table 8. It is noted that as $s$ increases, the standard error as well as the confidence intervals widths increase which means that predictive densities get more and more larger as the order statistics become larger.

\begin{center}
\begin{tabular}{llllllllll}
\multicolumn{10}{l}{Table 8: The values of point predictions and $95\%$ PIs for $y_s$.}\\\hline
\multicolumn{1}{c}{} & \multicolumn{1}{c}{} & \multicolumn{3}{c}{95\% PIs} & \multicolumn{1}{c}{} & \multicolumn{1}{c}{} & \multicolumn{3}{c}{95\% PIs} \\ \cline{3-5} \cline{8-10} \multicolumn{1}{c}{$s$} & \multicolumn{1}{c}{SEL} & \multicolumn{1}{c}{Lower} & \multicolumn{1}{c}{Upper} & \multicolumn{1}{c}{Length} & \multicolumn{1}{c}{$s$} & \multicolumn{1}{c}{SEL} & \multicolumn{1}{c}{Lower} & \multicolumn{1}{c}{Upper} & \multicolumn{1}{c}{Length} \\ \hline 1 & 0.182964 & 0.0049 & 0.64442 & 0.639522 & 12 & 2.3435 & 1.37295 & 3.4791 & 2.10615 \\ \hline 2 & 0.364917 & -0.04168 & 0.95305 & 0.994733 & 13 & 2.58715 & 1.56126 & 3.77906 & 2.2178 \\ \hline 3 & 0.546774 & 0.1234 & 1.22096 & 1.09755 & 14 & 2.84894 & 1.76337 & 3.56278 & 1.79941 \\ \hline 4 & 0.729419 & 0.21969 & 1.47136 & 1.25167 & 15 & 3.13419 & 1.98226 & 4.79029 & 2.80803 \\ \hline 5 & 0.913746 & 0.33041 & 1.71323 & 1.38282 & 16 & 3.45066 & 2.22227 & 3.64974 & 1.42747 \\ \hline 6 & 1.10068 & 0.45237 & 1.9516 & 1.49923 & 17 & 3.8104 & 2.48999 & 7.17473 & 4.68474 \\ \hline 7 & 1.29122 & 0.58385 & 2.18994 & 1.60609 & 18 & 4.23389 & 2.79605 & 7.36996 & 4.57391 \\ \hline 8 & 1.48646 & 0.72403 & 2.43111 & 1.70708 & 19 & 4.76065 & 3.15952 & 6.67321 & 3.51369 \\ \hline 9 & 1.68766 & 0.87265 & 2.67773 & 1.80509 & 20 & 5.48471 & 3.62082 & 7.7991 & 4.17828 \\ \hline 10 & 1.8963 & 1.02988 & 2.9325 & 1.90262 & 21 & 6.75234 & 4.30081 & 5.90783 & 1.60702 \\ \hline 11 & 2.11417 & 1.19631 & 3.19843 & 2.00212 & & & & & \\ \hline
\end{tabular}
\end{center}

\section{Simulation Study}

In this section, a simulation is performed to compare the estimators of the parameters on the basis of the samples generated from the LFR distribution. According to the Mean Square Errors (MSEs), average confidence interval lengths (ACILs) and the coverage probabilities (CPs), the performance of the competitive estimators have been compared.

Accordingly, a random sample of size $n$ has been generated from the LFR distribution for fixed values of its parameters $\alpha$ and $\beta$. To study the effect of the hybrid censoring scheme on the performance of the estimators, different combinations of the censoring parameters $(R,T)$ have been considered and simulation results are then summarized in Table 9 and Table 10.

\begin{center}
\begin{tabular}[c]{lcclllllclll}
\multicolumn{12}{l}{Table 9. MSE of MLEs and Bayesian estimates for the parameter $\alpha, \beta$ with $\alpha_{\circ
}=2, \beta_{\circ}=5$.}\\\cline{2-12}\cline{0-0}
&  &  &  &  & \multicolumn{3}{c}{$\alpha$} &  & \multicolumn{3}{c}{$\beta$%
}\\\cline{2-12}\cline{0-0}%
$n$ & $T$ & $R_{i}$ & $a_{r}$ &  & \ \ MLE &  & \ \ \ \ SEL &
\multicolumn{1}{l}{} & \ \ MLE &  & \ \ \ \ SEL\\\hline
$30$ & $3$ & $R_{2}=2,$ $R_{5}=1$ & $10$ &  &
\begin{tabular}
[c]{c}%
$2.3855$\\
$\left(  0.4031\right)  $%
\end{tabular}
&  &
\begin{tabular}
[c]{c}%
$2.3855$\\
$\left(  0.4031\right)  $%
\end{tabular}
& \multicolumn{1}{l}{} &
\begin{tabular}
[c]{c}%
$5.3006$\\
$\left(  0.8326\right)  $%
\end{tabular}
&  &
\begin{tabular}
[c]{c}%
$5.3007$\\
$\left(  0.8323\right)  $%
\end{tabular}
\\
&  & $R_{i}=0,$ $i\neq2,5$ & $18$ &  &
\begin{tabular}
[c]{c}%
$2.3491$\\
$\left(  0.3897\right)  $%
\end{tabular}
&  &
\begin{tabular}
[c]{c}%
$2.3491$\\
$\left(  0.3898\right)  $%
\end{tabular}
& \multicolumn{1}{l}{} &
\begin{tabular}
[c]{c}%
$5.3329$\\
$\left(  0.803\right)  $%
\end{tabular}
&  &
\begin{tabular}
[c]{c}%
$5.3329$\\
$\left(  0.8032\right)  $%
\end{tabular}
\\
&  &  & $25$ &  &
\begin{tabular}
[c]{c}%
$2.3789$\\
$\left(  0.4183\right)  $%
\end{tabular}
&  &
\begin{tabular}
[c]{c}%
$2.379$\\
$\left(  0.4183\right)  $%
\end{tabular}
& \multicolumn{1}{l}{} &
\begin{tabular}
[c]{c}%
$5.2453$\\
$\left(  0.7485\right)  $%
\end{tabular}
&  &
\begin{tabular}
[c]{c}%
$5.2453$\\
$\left(  0.7484\right)  $%
\end{tabular}
\\
& $5$ & $R_{2}=2,$ $R_{5}=1$ & $10$ &  &
\begin{tabular}
[c]{c}%
$2.3743$\\
$\left(  0.4287\right)  $%
\end{tabular}
&  &
\begin{tabular}
[c]{c}%
$2.3743$\\
$\left(  0.4287\right)  $%
\end{tabular}
& \multicolumn{1}{l}{} &
\begin{tabular}
[c]{c}%
$5.2584$\\
$\left(  0.7683\right)  $%
\end{tabular}
&  &
\begin{tabular}
[c]{c}%
$5.2583$\\
$\left(  0.7685\right)  $%
\end{tabular}
\\
&  & $R_{i}=0,$ $i\neq2,5$ & $18$ &  &
\begin{tabular}
[c]{c}%
$2.3768$\\
$\left(  0.4211\right)  $%
\end{tabular}
&  &
\begin{tabular}
[c]{c}%
$2.3768$\\
$\left(  0.4211\right)  $%
\end{tabular}
& \multicolumn{1}{l}{} &
\begin{tabular}
[c]{c}%
$5.2689$\\
$\left(  0.7412\right)  $%
\end{tabular}
&  &
\begin{tabular}
[c]{c}%
$5.269$\\
$\left(  0.7412\right)  $%
\end{tabular}
\\
&  &  & $25$ &  &
\begin{tabular}
[c]{c}%
$2.3799$\\
$\left(  0.3873\right)  $%
\end{tabular}
&  &
\begin{tabular}
[c]{c}%
$2.3798$\\
$\left(  0.3873\right)  $%
\end{tabular}
& \multicolumn{1}{l}{} &
\begin{tabular}
[c]{c}%
$5.3006$\\
$\left(  0.8011\right)  $%
\end{tabular}
&  &
\begin{tabular}
[c]{c}%
$5.3005$\\
$\left(  0.8011\right)  $%
\end{tabular}
\\
$40$ & $3$ & $R_{4}=2,$ $R_{9}=1$ & $18$ &  &
\begin{tabular}
[c]{c}%
$2.2745$\\
$\left(  0.2968\right)  $%
\end{tabular}
&  &
\begin{tabular}
[c]{c}%
$2.2745$\\
$\left(  0.2968\right)  $%
\end{tabular}
& \multicolumn{1}{l}{} &
\begin{tabular}
[c]{c}%
$5.2792$\\
$\left(  0.7653\right)  $%
\end{tabular}
&  &
\begin{tabular}
[c]{c}%
$5.2792$\\
$\left(  0.7656\right)  $%
\end{tabular}
\\
&  & $R_{i}=0,$ $i\neq2,9$ & $30$ &  &
\begin{tabular}
[c]{c}%
$2.2849$\\
$\left(  0.2944\right)  $%
\end{tabular}
&  &
\begin{tabular}
[c]{c}%
$2.2848$\\
$\left(  0.2944\right)  $%
\end{tabular}
& \multicolumn{1}{l}{} &
\begin{tabular}
[c]{c}%
$5.2671$\\
$\left(  0.7794\right)  $%
\end{tabular}
&  &
\begin{tabular}
[c]{c}%
$5.2673$\\
$\left(  0.7795\right)  $%
\end{tabular}
\\
&  &  & $35$ &  &
\begin{tabular}
[c]{c}%
$2.2443$\\
$\left(  0.2731\right)  $%
\end{tabular}
&  &
\begin{tabular}
[c]{c}%
$2.2443$\\
$\left(  0.2731\right)  $%
\end{tabular}
& \multicolumn{1}{l}{} &
\begin{tabular}
[c]{c}%
$5.329$\\
$\left(  0.8159\right)  $%
\end{tabular}
&  &
\begin{tabular}
[c]{c}%
$5.329$\\
$\left(  0.816\right)  $%
\end{tabular}
\\
& $5$ & $R_{4}=2,$ $R_{9}=1$ & $18$ &  &
\begin{tabular}
[c]{c}%
$2.2781$\\
$\left(  0.2949\right)  $%
\end{tabular}
&  &
\begin{tabular}
[c]{c}%
$2.2781$\\
$\left(  0.2949\right)  $%
\end{tabular}
& \multicolumn{1}{l}{} &
\begin{tabular}
[c]{c}%
$5.3002$\\
$\left(  0.8198\right)  $%
\end{tabular}
&  &
\begin{tabular}
[c]{c}%
$5.3001$\\
$\left(  0.8198\right)  $%
\end{tabular}
\\
&  & $R_{i}=0,$ $i\neq2,9$ & $30$ &  &
\begin{tabular}
[c]{c}%
$2.255$\\
$\left(  0.2914\right)  $%
\end{tabular}
&  &
\begin{tabular}
[c]{c}%
$2.255$\\
$\left(  0.2915\right)  $%
\end{tabular}
& \multicolumn{1}{l}{} &
\begin{tabular}
[c]{c}%
$5.3123$\\
$\left(  0.7931\right)  $%
\end{tabular}
&  &
\begin{tabular}
[c]{c}%
$5.3124$\\
$\left(  0.7931\right)  $%
\end{tabular}
\\
&  &  & $35$ &  &
\begin{tabular}
[c]{c}%
$2.2821$\\
$\left(  0.3118\right)  $%
\end{tabular}
&  &
\begin{tabular}
[c]{c}%
$2.2821$\\
$\left(  0.3118\right)  $%
\end{tabular}
& \multicolumn{1}{l}{} &
\begin{tabular}
[c]{c}%
$5.2793$\\
$\left(  0.7458\right)  $%
\end{tabular}
&  &
\begin{tabular}
[c]{c}%
$5.2793$\\
$\left(  0.7458\right)  $%
\end{tabular}
\\\hline
\end{tabular}

\begin{tabular}[c]{lcclllllclll}
\multicolumn{12}{l}{Table 10. ACIL and CP of $95\%$ CIs for the parameters $\alpha$ and $\beta$}\\\hline
&  &  &  &  & \multicolumn{3}{c}{$\alpha$} &  & \multicolumn{3}{c}{$\beta$%
}\\\cline{2-12}\cline{0-0}
&  &  &  &  & MLE &  & MCMC & \multicolumn{1}{l}{} & MLE &  &
MCMC\\\cline{2-12}\cline{0-0}%
$n$ & $T$ & $R_{i}$ & $a_{r}$ &  &  &  &  & \multicolumn{1}{l}{} &  &  &
\\\hline
$30$ & $3$ & $R_{2}=2,$ $R_{5}=1$ & $10$ &  &
\begin{tabular}
[c]{c}%
$3.397$\\
$\left(  0.9456\right)  $%
\end{tabular}
&  &
\begin{tabular}
[c]{c}%
$0.0015$\\
$\left(  0.9358\right)  $%
\end{tabular}
& \multicolumn{1}{l}{} &
\begin{tabular}
[c]{c}%
$13.9187$\\
$\left(  0.9312\right)  $%
\end{tabular}
&  &
\begin{tabular}
[c]{c}%
$0.0061$\\
$\left(  0.9294\right)  $%
\end{tabular}
\\
&  & $R_{i}=0,$ $i\neq2,5$ & $18$ &  &
\begin{tabular}
[c]{c}%
$3.3796$\\
$\left(  0.9433\right)  $%
\end{tabular}
&  &
\begin{tabular}
[c]{c}%
$0.0014$\\
$\left(  0.9329\right)  $%
\end{tabular}
& \multicolumn{1}{l}{} &
\begin{tabular}
[c]{c}%
$13.8409$\\
$\left(  0.971\right)  $%
\end{tabular}
&  &
\begin{tabular}
[c]{c}%
$0.0058$\\
$\left(  0.9606\right)  $%
\end{tabular}
\\
&  &  & $25$ &  &
\begin{tabular}
[c]{c}%
$3.3822$\\
$\left(  0.9738\right)  $%
\end{tabular}
&  &
\begin{tabular}
[c]{c}%
$0.0015$\\
$\left(  0.9666\right)  $%
\end{tabular}
& \multicolumn{1}{l}{} &
\begin{tabular}
[c]{c}%
$13.8007$\\
$\left(  0.9538\right)  $%
\end{tabular}
&  &
\begin{tabular}
[c]{c}%
$0.0059$\\
$\left(  0.9483\right)  $%
\end{tabular}
\\
& $5$ & $R_{2}=2,$ $R_{5}=1$ & $10$ &  &
\begin{tabular}
[c]{c}%
$3.3913$\\
$\left(  0.9329\right)  $%
\end{tabular}
&  &
\begin{tabular}
[c]{c}%
$0.0015$\\
$\left(  0.9316\right)  $%
\end{tabular}
& \multicolumn{1}{l}{} &
\begin{tabular}
[c]{c}%
$13.8519$\\
$\left(  0.9489\right)  $%
\end{tabular}
&  &
\begin{tabular}
[c]{c}%
$0.006$\\
$\left(  0.9508\right)  $%
\end{tabular}
\\
&  & $R_{i}=0,$ $i\neq2,5$ & $18$ &  &
\begin{tabular}
[c]{c}%
$3.3932$\\
$\left(  0.9377\right)  $%
\end{tabular}
&  &
\begin{tabular}
[c]{c}%
$0.0014$\\
$\left(  0.9713\right)  $%
\end{tabular}
& \multicolumn{1}{l}{} &
\begin{tabular}
[c]{c}%
$13.8733$\\
$\left(  0.9568\right)  $%
\end{tabular}
&  &
\begin{tabular}
[c]{c}%
$0.006$\\
$\left(  0.9436\right)  $%
\end{tabular}
\\
&  &  & $25$ &  &
\begin{tabular}
[c]{c}%
$3.3956$\\
$\left(  0.9400\right)  $%
\end{tabular}
&  &
\begin{tabular}
[c]{c}%
$0.0015$\\
$\left(  0.9406\right)  $%
\end{tabular}
& \multicolumn{1}{l}{} &
\begin{tabular}
[c]{c}%
$13.9008$\\
$\left(  0.9576\right)  $%
\end{tabular}
&  &
\begin{tabular}
[c]{c}%
$0.0059$\\
$\left(  0.9695\right)  $%
\end{tabular}
\\
$40$ & $3$ & $R_{4}=2,$ $R_{9}=1$ & $18$ &  &
\begin{tabular}
[c]{c}%
$2.9366$\\
$\left(  0.9303\right)  $%
\end{tabular}
&  &
\begin{tabular}
[c]{c}%
$0.0013$\\
$\left(  0.9571\right)  $%
\end{tabular}
& \multicolumn{1}{l}{} &
\begin{tabular}
[c]{c}%
$11.8872$\\
$\left(  0.9504\right)  $%
\end{tabular}
&  &
\begin{tabular}
[c]{c}%
$0.005$\\
$\left(  0.9483\right)  $%
\end{tabular}
\\
&  & $R_{i}=0,$ $i\neq2,9$ & $30$ &  &
\begin{tabular}
[c]{c}%
$2.9514$\\
$\left(  0.9273\right)  $%
\end{tabular}
&  &
\begin{tabular}
[c]{c}%
$0.0013$\\
$\left(  0.9513\right)  $%
\end{tabular}
& \multicolumn{1}{l}{} &
\begin{tabular}
[c]{c}%
$11.9389$\\
$\left(  0.9533\right)  $%
\end{tabular}
&  &
\begin{tabular}
[c]{c}%
$0.0051$\\
$\left(  0.9269\right)  $%
\end{tabular}
\\
&  &  & $35$ &  &
\begin{tabular}
[c]{c}%
$2.929$\\
$\left(  0.9694\right)  $%
\end{tabular}
&  &
\begin{tabular}
[c]{c}%
$0.0012$\\
$\left(  0.9742\right)  $%
\end{tabular}
& \multicolumn{1}{l}{} &
\begin{tabular}
[c]{c}%
$11.8674$\\
$\left(  0.9462\right)  $%
\end{tabular}
&  &
\begin{tabular}
[c]{c}%
$0.005$\\
$\left(  0.9666\right)  $%
\end{tabular}
\\
& $5$ & $R_{4}=2,$ $R_{9}=1$ & $18$ &  &
\begin{tabular}
[c]{c}%
$2.9452$\\
$\left(  0.9688\right)  $%
\end{tabular}
&  &
\begin{tabular}
[c]{c}%
$0.0013$\\
$\left(  0.9746\right)  $%
\end{tabular}
& \multicolumn{1}{l}{} &
\begin{tabular}
[c]{c}%
$11.9416$\\
$\left(  0.9648\right)  $%
\end{tabular}
&  &
\begin{tabular}
[c]{c}%
$0.0051$\\
$\left(  0.9324\right)  $%
\end{tabular}
\\
&  & $R_{i}=0,$ $i\neq2,9$ & $30$ &  &
\begin{tabular}
[c]{c}%
$2.9269$\\
$\left(  0.9690\right)  $%
\end{tabular}
&  &
\begin{tabular}
[c]{c}%
$0.0013$\\
$\left(  0.9637\right)  $%
\end{tabular}
& \multicolumn{1}{l}{} &
\begin{tabular}
[c]{c}%
$11.8574$\\
$\left(  0.9538\right)  $%
\end{tabular}
&  &
\begin{tabular}
[c]{c}%
$0.005$\\
$\left(  0.9512\right)  $%
\end{tabular}
\\
&  &  & $35$ &  &
\begin{tabular}
[c]{c}%
$2.9454$\\
$\left(  0.9732\right)  $%
\end{tabular}
&  &
\begin{tabular}
[c]{c}%
$0.0013$\\
$\left(  0.9461\right)  $%
\end{tabular}
& \multicolumn{1}{l}{} &
\begin{tabular}
[c]{c}%
$11.9353$\\
$\left(  0.9467\right)  $
\end{tabular}
&  &
\begin{tabular}
[c]{c}%
$0.0051$\\
$\left(  0.9373\right)  $
\end{tabular}
\\\hline
\end{tabular}
\ \
\end{center}

From the reported values, we observed the following:

\begin{itemize}
  \item As the sample size increases, the MSEs for most of the estimators decreases, for the same values of $\alpha$ and $\beta$.
  \item For different values of $n,~T~\text{and}~a_r$, the MSEs of Bayes estimators and MLEs are nearly the same.
  \item The Bayes estimators obtained are nearly the same to the MLEs. Also it is noticed that providing prior information gives more accurate results for the Bayesian estimates of the parameters.
  \item The credible intervals width is less than that of the asymptotic confidence intervals in all the cases. However, it is noted that the confidence intervals width decreases as the sample size increases.
\end{itemize}

\section{Conclusion}
The Linear failure rate distribution is well-fitted to some real-life data which leads to obtaining good estimators for its parameters. Multiply Type-II hybrid censoring scheme is a good description model for experiments that lack the exact registration of the failure times for some units. A prediction study is applied and proven that the predictive sample simulates the observed sample. It is recommended to apply the Multiply hybrid censoring scheme to clinical data especially in cases of pandemics when it is difficult for researchers to register all the failure times exactly.


\begin{thebibliography}{99}

\bibitem{1}
E. A. Ahmed, \emph{Bayesian estimation based on progressive Type-II censoring from two-parameter bathtub-shaped lifetime model: a Markov chain Monte Carlo approach}, Journal of Applied Statistics 41 (4) 2014, pp. 752—768.

\bibitem{2}
E. A. Ahmed, \emph{Estimation and prediction for the generalized inverted exponential distribution based on progressively first-failure-censored data with application}, Journal of Applied Statistics 2014. Available at http://dx.doi.org/10.1080/02664763.2016.1214692.

\bibitem{3}
B. Al-zahrani and M. Gindwan, \emph{Parameter estimation of a two-parameter Lindley distribution under hybrid censoring}, International Journal of System Assurance Engineering and Management 5 (2014), pp. 628--636.

\bibitem{4}
A. Azzalini, \emph{Statistical inference based on the likelihood}, Chapman and Hall, London (1996).

\bibitem{5}
A. Banerjee and D. Kundu, \emph{Inference based on Type-II hybrid censored data from a Weibull distribution}, IEEE Transactions on Reliability 57 (2008), pp. 369--378.

\bibitem{6}
S. M. Chen and G. K. Bhattacharyya, \emph{Exact confidence bound for an exponential under hybrid censoring}, Communication in Statistics - Theory and Methods 17 (1988), pp. 1858--1870.

\bibitem{7}
Chien-Wei Wu, Ming-Hung Shu and Yu-Ning Chang, \emph{Variable-sampling plans based on lifetime-performance index under exponential distribution with censoring and its extensions}, Applied Mathematical Modelling 55 (2018), pp. 81--93.

\bibitem{8}
A. Childs, B. Chandrasekhar, N. Balakrishnan and D. Kundu, \emph{Exact likelihood inference based on type-I and type-II hybrid censored samples from the exponential distribution}, Annals of the Institute of Statistical Mathematics 55 (2003), pp. 319--330.

\bibitem{9}
Chin-Diew Lai and Min Xie, \emph{Stochastic Ageing and Dependence for Reliability}, Springer, New York, NY.

\bibitem{10}
M. Y. Danish, I. A. Arshad and M. Aslam, \emph{Bayesian inference for the randomly censored Burr-type XII distribution}, Journal of Applied Statistics 45 (2) 2018, pp. 270—283.

\bibitem{11}
H.A. David and H.N. Nagaraja, \emph{Order Statistics}, 3rd ed., John Wiley \& Sons, New York, 2003.

\bibitem{12}
S. Dey and B. Pradhan, \emph{Generalized inverted exponential distribution under hybrid censoring}, Statistical Methodology 18 (2014), pp. 101--114.

\bibitem{13}
S. Dubey, B. Pradhan and D. Kundu (2011), \emph{Parameter estimation of the hybrid censored log-normal distribution}, Journal of Statistical Computation and Simulation 81 (2011), pp. 275--287.

\bibitem{14}
R. M. EL-Sagheer, \emph{Estimation of parameters of Weibull-Gamma distribution based on progressively censored data}, Statistical Papers 59 (2018), pp. 725—757. Available at https://doi.org/10.1007/s00362-016-0787-2.

\bibitem{15}
A. Ganguly, S. Mitra, D. Samanta and D. Kundu, \emph{Exact inference for the two parameter exponential distribution under Type-II hybrid censoring}, Journal of Statistical Planning and Inference 142 (2012), pp. 613--625.

\bibitem{16}
P. K. Gupta and B. Singh, \emph{Parameter estimation of Lindley distribution with hybrid censored data}, International Journal of System Assurance Engineering and Management, 4 (2013), pp. 378--385.

\bibitem{17}
Y. E. Jeon, S. B. Kang, \emph{Bayesian estimation for the exponential distribution based on generalized multiply Type-II hybrid censoring}, Communications for Statistical Applications and Methods 27 (4) 2020, pp. 413—430.

\bibitem{18}
Y. E. Jeon, S. B. Kang, \emph{Estimation for the half-logistic distribution based on multiply Type-II hybrid censoring}, Physica A 550 (2020). Available at https://doi.org/10.1016/j.physa.2020.124501.

\bibitem{19}
N. L. Johnson and S. Kotz, \emph{Encyclopedia of Statistical Sciences}, Wiley, New York, (1983).

\bibitem{20}
J. P. Keating, R. E. Glaser and N. S. Ketchum, \emph{Testing hypotheses about the shape of a gamma distribution}, Technometrics 32 (1990), 67-82.

\bibitem{21}
A. Kohansal, S. Rezakhah and E. Khorram, \emph{Parameter estimation of Type-II hybrid censored weighted exponential distribution}, Communications in Statistics - Simulation and Computation 44 (2015), pp. 1273--1299.

\bibitem{22}
D. Kundu, \emph{On hybrid censored Weibull distribution}, Journal of Statistical Planning and Inference 137 (2007), pp. 2127--2142.

\bibitem{23}
J. F. Lawless, \emph{Statistical models and methods for lifetime data}, John Wiley and Sons, New York (1982).

\bibitem{24}
WH. Lee and KJ. Lee, \emph{Estimating the parameter of an exponential distribution under multiply Type-II hybrid censoring}, Journal of the Korean Data and Information Science Society 29 (2018), pp. 807--814.

\bibitem{25}
M. A. W. Mahmoud , D. A. Ramadan and M. M. M. Mansour, \emph{Estimation of lifetime parameters of the modified extended exponential distribution with application to a mechanical model}, Communications in Statistics - Simulation and Computation 2020. Available at https://doi.org/10.1080/03610918.2020.1821887.

\bibitem{26}
M. M. M. Mansour and D. A. Ramadan, \emph{STATISTICAL INFERENCE TO THE PARAMETERS OF THE MODIFIED EXTENDED EXPONENTIAL DISTRIBUTION UNDER THE TYPE-II HYBRID CENSORING SCHEME}, Journal of Applied Probability and Statistics 15 (2) 2020, pp. 19—44.

\bibitem{27}
S. Nadarajaha and S. A. A. Bakar, \emph{An exponentiated geometric distribution}, Applied Mathematical Modelling 40 (2016), pp. 6775--6784.

\bibitem{28}
S. Park and N. Balakrishnan, \emph{A very flexible hybrid censoring scheme and its Fisher information}, Journal of Statistical Computation and Simulation 82 (2012), pp. 41--50.

\bibitem{29}
M. K. Rastogi and Y. M. Tripathi, \emph{Inference on unknown parameters of a Burr distribution under hybrid censoring}, Statistical Papers 54 (2013), pp. 619--643.

\bibitem{30}
R. Royall, \emph{Statistical Evidence: A Likelihood Paradigm}, Chapman and Hall, London (1997).

\bibitem{31}
Tian-Qun Xu and Yue-Peng Chen, \emph{Two-sided M-Bayesian credible limits of reliability parameters in the case of zero-failure data for exponential distribution}, Applied Mathematical Modelling 38 (2014), pp. 2586--2600.

\end{thebibliography}
\end{document}